\documentclass[12pt]{article}
\usepackage{amsmath}
\usepackage{amssymb}
\usepackage[mathscr]{eucal}
\usepackage{graphicx}
\usepackage{epstopdf} 
\usepackage{color}
\usepackage{xcolor}
\usepackage{float}

\usepackage{textcomp}
\usepackage{subfig}
\usepackage{placeins}
\usepackage{url}
\usepackage[normalem]{ulem}

\def\real{\mathbb{R}}
\def\nat{\mathbb{N}}

\newcommand{\T}{\mathcal{T}}

\newcommand{\skalar}[1]{\left\langle #1 \right\rangle}
\newcommand{\skalareh}[2]{\left\langle \{\nabla#1\},[#2] \right\rangle_{\mathcal E_h}}
\newcommand{\skalarehzero}[2]{\left\langle \{\nabla#1\},[#2] \right\rangle_{\mathcal E_h^0}}
\newcommand{\sskalareh}[2]{\left\langle [#1],[#2] \right\rangle_{\mathcal E_h}}
\newcommand{\sskalarehzero}[2]{\left\langle [#1],[#2] \right\rangle_{\mathcal E_h^0}}
\newcommand{\sskalarehpartial}[2]{\left\langle [#1],[#2] \right\rangle_{\mathcal E_h^\partial}}

\newcommand{\Th}{\mathcal T_h}

\allowdisplaybreaks
\newtheorem{Theorem}{\sc Theorem}

\newtheorem{Proposition}[Theorem]{\sc Proposition}
\newtheorem{Lemma}[Theorem]{\sc Lemma}
\newtheorem{Corollary}[Theorem]{\sc Corollary}

\def\sqr#1#2{{
    \vcenter{
         \vbox{\hrule height.#2pt
               \hbox{\vrule width.#2pt height#1pt \kern#1pt
                     \vrule width.#2pt
               }
               \hrule height.#2pt
         }
    }
}}

\def\bar{\overline}
\def\real{\mathbb{R}}
\def\nat{\mathbb{N}}

\newcommand{\R}{{\if mm {\rm I}\mkern -3mu{\rm R}\else \leavevmode
\hbox{I}\kern -.17em\hbox{R} \fi}}

\def\lista#1
{{ \itemindent 0.0cm \labelsep .2cm \leftmargin 0.8cm \rightmargin
0.0cm \labelwidth 0.6cm \topsep 0.0mm
\parsep 0.0mm
\itemsep 0.0mm
\begin{list}{}
{ \setlength{\leftmargin}{1.0cm} \setlength{\rightmargin}{0.0cm}
\setlength{\parsep}{.0mm} \setlength{\topsep}{.0mm}
\setlength{\parskip}{.0cm} \setlength{\itemsep}{.0cm} }
{#1}\end{list}} }

\begin{document}

\title{ Discontinuous Virtual Element Method for an elliptic variational inequality of the second kind
		 \ }

\author{
Krzysztof Bartosz$^{\,1}$,  \
Pawe{\l} Szafraniec$^{\,1}$ \ 
\\  \\
{\small $^1$ Jagiellonian University, Faculty of Mathematics and Computer Science} \\
{\small ul. Lojasiewicza 6, 30348 Krakow, Poland} \\
}

\date{}
\maketitle

\vskip 4mm

\noindent {\bf Abstract.} In this work, we analyse a simplified frictional contact problem and its variational formulation that has a form of the  elliptic variational inequality of the second kind. For this problem, we consider a numerical approximation based on virtual element method (VEM). The novelty of our approach consists in the fact that we deal with a discontinuous version of the VEM. Our analysis concerns the error estimation between the exact solution of the simplified frictional contact problem and the approximate one.   

\vskip 4mm


\vskip 4mm


\thispagestyle{empty}

\





\section{Introduction}\label{Introduction}
\begin{figure}[H]
	\centering
	\includegraphics[width=14cm]{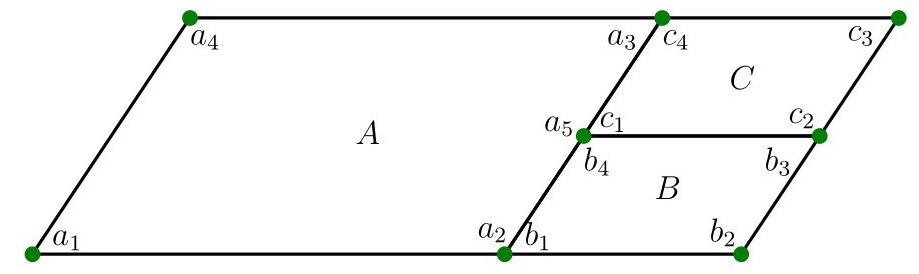}
	\caption{Hanging node $a_5=b_4=c_1$}\label{fig_1}
\end{figure}

The goal of this paper is to study an innovative numerical approach for solving an  elliptic variational inequality of the second kind. Our idea is based on a mixture of two strategies that represent a generalization of Finite Element Method (FEM). Here we briefly describe a motivation to use this concept. For solving elliptic problems, the FEM seems to be a natural choice, but in some situations it may be not effective. For example, it fails to handle hanging nodes that may occur during the mesh refinement (as depicted in Figure \ref{fig_1}). There are at least two ways to fix this issue. The first one is to apply the Discontinuous Galerkin Method (DGM). In this approach the approximate solution of a problem, related to an exemplary geometry presented in Figure \ref{fig_1}, is composed of three independent parts described on elements $A$, $B$ and $C$, respectively. Each of parts belongs to $Q_1(K)$, $K\in\{A,B,C\}$, where $Q_1(K)$ denotes the space of first order rectangular shape functions (see  \cite[Chapter 3.5]{BRESCO}) modified up to some affine transformation (see \cite[page 122]{50_lines} for details). For reference on the DGM, see e.g. \cite{Tutorial_DG}.
An alternative way of handling hanging nodes, is to apply the Virtual Element Method (VEM). In this approach, the element $A$ presented in Figure \ref{fig_1}, is treated as a pentagon with the vertices $a_1, a_2, a_5, a_3, a_4$. The solution still consists of three parts, and each of them belongs to the local virtual element space $V_1(K)$, $K\in\{A,B,C\}$, respectively (see e.g. \cite{Principles} for details). The degrees of freedom are represented by values of the function at each of nodes. The solution is continuous, as the global degrees of freedom corresponding to the nodes $a_5$, $b_4$ and $c_1$ coincide. This method introduced in a seminal paper \cite{Principles} permits the use of meshes with highly irregular shaped elements and arbitrary number of nodes. It became quickly very popular topic of research, both from practical and theoretical point of view. We recommend \cite{Hitchhiker} to start with the VEM, for a quick tutorial on how to implement the VEM see e.g. \cite{VEM-50}. The VEM has been implemented to various problems in elasticity and other physical and engineering problems, see e.g. \cite{VEM_elasticity, Principles_mixed, Conforming_and_non, FWH, VEM_elliptic_HVI, Wriggers, WWH}.

 However, since the publication of \cite{Principles} some small, yet important gaps in the theoretical foundations have been found. The main issue was the lack of affine equivalence between elements, standard for the FEM, but in the case of arbitrary polygons becomes difficult. In the aforementioned paper \cite{Principles}, authors prove the result without this assumption, and resolve the problem, referring to \cite{BRESCO}. But the careful analysis of \cite{BRESCO} ensures that the affine equivalence of elements is required there. Roughly speaking, two elements are affine equivalent if there exists a nonsingular affine transformation that maps one element to another such that the degrees of freedom associated to the first element correspond to those of the second one via the composition with the aforementioned affine transformation (see  \cite[Definition 3.4.1]{BRESCO} for details). 
 
 Coming back to the Figure \ref{fig_1}, its clear, that the virtual element $A=(a_1,a_2,a_5,a_3,a_4)$ is not affine equivalent with the others, as it simply has more degrees of freedom. Hence, even if there exists an affine transformation that maps the parallelogram $A$ to  parallelogram $B$ or $C$, it cannot preserve the required correspondence between the degrees of freedom. To fix this issue, we simply remove the node $a_5$ and write  $A=(a_1,a_2,a_3,a_4)$. Now all virtual elements are affine equivalent and one can take all advantages of this fact. But the price, one has to pay for this simplification, is a discontinuity that occurs along the common edge. Hence, we are led to the Discontinuous Virtual Element Method (DVEM). This idea have been started in \cite{DGVG}, but to our knowledge has not been pursued further.\\
 We mention that in specific cases the DVEM reduces to the aforementioned DGM. It happens if each element $K$ is a rectangle parallel to the axes of the system of coordinates. Then one has $V_1(K)=Q_1(K)$ and both methods coincide. However, in general case, the DVEM is a significant novelty in comparison to the DGM. 
  
This paper concentrates on combining both the VEM and the DG methods for a second kind variational inequality describing simplified friction problem. This problem was analysed by the DGM in \cite{F_Wang_DG_Ell_VI} and by the VEM in \cite{Fei_Wang}, hence, we generalize these two results by mixing both ideas. 
 Our approach is the following: we assume explicitly affine equivalence of elements, which restricts us to using only one fixed type of polygon on the whole mesh. It allows to base on most classical results for example \cite{BRESCO}. On the down side, it may produce hanging nodes, but this is the point where discontinuous part comes in, allowing us to bypass this issue. In our opinion, while slightly more technical to some extent, this approach could be of value, as it is more similar to the classical FEM and hence more intuitive. Finally, we mention that using some more sophisticated argumentation  it may be possible to obtain our result without the restrictive assumption on affine equivalence. It can be provided by the so called virtual triangulation of elements (cf. Assumption A1 in \cite{ChenHuang}). Nevertheless, for the sake of simplicity and consistency of the proof, we decided to base it on the affine equivalence property.

The structure of the paper is as follows. In Section 2 we recall the concept of variational inequalities (VI) and present the classical elliptic problem that leads to the VI. In Section 3 we introduce suitable discretization of the domain and present preliminary material concerning the DVEM. In the last section we carefully prove all inequalities needed for the main theorem of the paper, namely Theorem 9, which gives a linear bound for the difference between approximate and true solution. In our view this approach should be pursued further. Introducing for example nonconvex potential through hemivariational inequalities  see e.g. \cite{MS} for reference, may be of great interest.

\section{Problem formulation}
First, we introduce the notion of a variational inequality. Let $V$ be a real Hilbert space with scalar product $(\cdot,\cdot)$ and associated norm $\|\cdot\|$ and $V^*$ be its dual. Consider the following problem: find $u\in V$ such that
\begin{align}\label{EVI}
a(u,v-u)+j(v)-j(u)\geq L(u-v)\quad \text{for\, all}\,\,\, v\in V,
\end{align} 
where $a\colon V\times V\to\real$, $L\colon V\to\real$ and $j\colon V\to\bar{\real}=\real\cup\{\infty\} $. Problem (\ref{EVI}) is known as an elliptic variational inequality (EVI) of the second kind. We recall a result on existence and uniqueness of its solution (cf. \cite[Chapter 1, Theorem 4.1]{Glo}). 
\begin{Theorem}\label{Theorem_0}
	Assume that
	\begin{itemize}
		\item $a\colon V\times V\to\real$ is a bilinear, continuous and $V$-elliptic form,
		\item $L\colon V\to R$ is a linear and continuous functional (i.e. $f\in V^*$),
		\item $j\colon V\to\bar{\real}=\real\cup\{\infty\}$ is a convex, lower semicontinuous (l.s.c) and proper functional. 
	\end{itemize}
Then problem (\ref{EVI}) has a unique solution. 
\end{Theorem} 

Now we describe the physical problem that can be modelled by means of the EVI. We consider a bounded domain $\Omega\subset\real^2$, whose regular boundary $\partial\Omega=\bar{\Gamma}_1\cup\bar{\Gamma}_2$ is a union of two disjoint parts $\Gamma_1$ and $\Gamma_2$. Let $n$ denote the unit outward normal vector on the boundary $\partial \Omega$. Let $g>0$ and $f\in L^2(\Omega)$ be given.  We consider the following problem.\\[3mm] 
\noindent {\bf Problem ($P$).} {\it Find $u\in H^1(\Omega)$ such that there exists $\lambda\in L^{\infty}(\Gamma_2)$ satisfying
\begin{align}
-\Delta u+u=f\qquad & \text{in}\,\,\,\,\,\,\Omega,\label{1}\\[1mm]
\nabla u\cdot n=0\qquad & \text{on}\,\,\,\,\,\,\Gamma_1,\label{2}\\[1mm]
\nabla u\cdot n+g\lambda=0 \qquad & \text{on}\,\,\,\,\,\,\Gamma_2,\label{3}\\[1mm]
|\lambda|\leq 1,\,\, \lambda u=|u|\qquad & \text{a.e. on} \,\,\,\,\,\, \Gamma_2.\label{3b}
\end{align}
} 

\noindent Problem $(P)$ is referred to as a simplified frictional contact problem and it was studied in \cite{Glo} for instance. \\
Now we deal with variational formulation of Problem ($P$). Let $V=H^1(\Omega)$ and let
\begin{align*}
&a(u,v)=\int_{\Omega}\nabla u\cdot \nabla v\,dx+\int_{\Omega}u\cdot v\,dx\qquad \mbox{for all} ~ \,\, u,v\in V,\\[2mm]
&(f,v)=\int_{\Omega}fv\,dx\qquad \mbox{for all} ~ \,\,v\in V,\\[2mm]
&j(v)=\int_{\Gamma_2}g|v|\,ds\qquad \mbox{for all} ~ \,\, v\in V.
\end{align*}
We assume that the space $V$ is endowed with the classical scalar product and its associated norm, i.e. $(u,v)_V=a(u,v)$ and $\|u\|^2_V=a(u,u)$ for all $u,v\in V$.\\ 
 
\noindent The variational formulation of Problem ($P$) reads as follows.\\
 
\noindent {\bf Problem ($V$).} {\it Find $u\in V$ such that
	\begin{align}
	a(u,v-u)+j(v)-j(u)\geq (f,v-u)\qquad \mbox{for all} ~ \,\, v\in V.\label{eq_1}
	\end{align}
}\noindent
Inequality (\ref{eq_1}) can be obtained from (\ref{1}) by multiplying it by $v-u$, integrating over $\Omega$, applying Green formula and boundary conditions (\ref{2})-(\ref{3}) and using properties (\ref{3b}). Hence, every solution of Problem ($P$) is a solution of Problem $(V)$. The inverse relation is not trivial, however, using analogous arguments as in the proof of  \cite[Chapter 2, Theorem 5.3]{Glo} and considering   \cite[Chapter 2, Remark 5.4]{Glo}, we conclude, that every solution $u$ of Problem ($V$) is also a solution of Problem ($P$). Moreover,  its corresponding multiplier $\lambda$ is unique as well. Furthermore (see \cite[Chapter 2, Theorem 5.2]{Glo}), every solution $u$ of Problem ($V$) has regularity $u\in H^2(\Omega)$. Moreover, it is easy to check, that in our case, all assumptions of Theorem \ref{Theorem_0} are satisfied. Hence, it follows that Problem $(V)$ has a unique solution. Summarizing, we have.
\begin{Corollary}\label{Corollary_1}
	Problem $(V)$ has a unique solution $u\in H^2(\Omega)$. Moreover $u$ is a unique solution of Problem ($P$) with a unique multiplier $\lambda$.  
\end{Corollary}

\section{Discretization}
In this section we define a suitable discretization of the domain $\Omega$. Next, we present preliminary material concerning the DGM and the VEM that will be useful in the numerical analysis of our problem.\\    	 	
Let $\{{\mathcal T_h}\}_h$ be a sequence of decompositions of $\Omega$ into elements $K$ that has the following properties.
\begin{itemize}
	\item $\bar{\Omega}=\bigcup_{K\in \mathcal T_h }K$. 
	\item Each element $K\in \mathcal T_h$ is a simple polygon.
	\item For distinct $K_1, K_2\in {\mathcal T_h}$,\, $\overset{\circ}{K_1}\cap\overset{\circ}{K_2}=\emptyset$.
\end{itemize}
	
We remark at this point that we do not assume other common hypothesis, namely: {\it For distinct $K_1, K_2\in {\mathcal T_h}$,\, ${K_1}\cap{K_2}$ is empty, or a common vertex, or a common side of $K_1$ and $K_2$}. This is due to the fact, that the main idea of the paper is to deal with situation where this may not be satisfied.

The symbol $h$ denotes the parameter of discretization defined by  
		\begin{align*}
h=\max_{K\in\T_h}h_K,\,\,\,\,\text{where}\,\,\,\, h_K=\max_{x,y\in K}\|x-y\|_{\real^2}. 
\end{align*}
We introduce the following assumptions.
\begin{itemize}
	\item [$(C_0)$] There is a reference element $\hat{K}$ such that for all $K\in \mathcal T_h$ there exists a bijective affine function $F_K\colon \hat{K}\to K$.
	\item[$(C_1)$] There exists a real number $\gamma_1>0$ such that each element $K\in \Th $ 
	is star-shaped with respect to a disk of radius $\rho_K\geq \gamma_1 h_K$.
	\item[$(C_2)$] There exists a real number $\gamma_2>0$ such that for each element $K\in \Th $, the
	distance between any two vertices of $K$ is $\geq \gamma_2 h_K$.
\end{itemize}
We use the following notation:
\begin{itemize}
	\item $\mathcal T_h$ - set of all elements,
	\item $\mathcal E_h$ - set of all edges,
	\item $\mathcal E_h^0$ - set of all internal edges,
	\item $\mathcal E_h^{\partial}$ - set of all boundary edges.
\end{itemize}
We define the space
\begin{align*}
H^m({\mathcal T_h})=\{ v\in L^2(\Omega):\,\,v_{|K}\in H^m(K)\,\,\,\,\mbox{for all} ~  K\in  {\mathcal T_h} \},\quad m\in\nat.
\end{align*}
For $v\in H^1({\mathcal T_h})$ and $\tau\in (H^1({\mathcal T_h}))^2$ we define the averages and jumps on the internal edge $e=K_1\cap K_2$ by te formulas
\begin{align*}
&\{v\}=\frac{v^1+v^2}{2},\qquad\qquad 
 [v]=v^1n^1+v^2n^2,\\[2mm]
&\{\tau\}=\frac{\tau^1+\tau^2}{2} ,\qquad\qquad [\tau]=\tau^1\cdot n^1+\tau^2\cdot n^2.
\end{align*}
Here $v^i, \tau^i$ denote the values of $v_{|K_i}$, $\tau_{|K_i}$ at the edge $e$ in the sense of trace operator for $i=1,2$. Moreover, $n^i$ denotes the unit outward normal vector to the element $K_i$ at the boundary $e$, $i=1,2$. Clearly, we have $n^1=-n^2$. On the boundary edge $e$, we define 
\begin{align*}
[v]=vn,\quad\{\tau\}=\tau.
\end{align*}
We also define for all $ v\in H^2({\mathcal T_h})$ and $w\in H^1({\mathcal T_h})$ 
\begin{align*}
\skalareh{v}{w}=\sum_{e\in{\mathcal E_h}}\int_e\{\nabla v\}[w]ds,\\
\skalarehzero{v}{w}=\sum_{e\in{\mathcal E_h^0}}\int_e\{\nabla v\}[w]ds,
\end{align*}
and for all $ v,w\in H^1({\mathcal T_h})$ we write
\begin{align*}
&\sskalareh{v}{w}=\sum_{e\in{\mathcal E_h}}\frac{1}{|e|}\int_e[v][w]ds,\\
&\sskalarehzero{v}{w}=\sum_{e\in{\mathcal E_h^0}}\frac{1}{|e|}\int_e[v][w]ds,\\
&\sskalarehpartial{v}{w}=\sum_{e\in{\mathcal E_h^\partial}}\frac{1}{|e|}\int_e[v][w]ds,\\
&\|[v]\|^2_{0,\partial K}=\sum_{e\in \partial K}\frac{1}{|e|}\int_e|[v]|^2ds\quad\text{for all}\,\,K\in{\mathcal T_h}.
\end{align*}

We introduce the notation 
\begin{align*}
a^K(u,v)=\int_{K}\nabla u\cdot\nabla v+uv\,dx \quad\text{for all}\,\,u,v\in H^1(K),
\end{align*}
for all $K\in {\mathcal T_h}$. Hence, we can extend the notion of bilinear form $a$ as follows. Let $\tilde{a}\colon  H^1({\mathcal T_h})\times  H^1({\mathcal T_h})\to\real $ be defined by
\begin{align*}
\tilde{a}(u,v)=\sum_{K\in\Th}a^K(u,v)\quad\text{for all}\,\,u,v\in  H^1({\mathcal T_h}).
\end{align*}  
Clearly, we have $a(u,v)=\tilde{a}(u,v)$ for $u,v\in H^1(\Omega)$.

\begin{Lemma}\label{Lemma_1}
	If $u\in V$ is a solution of Problem $(V)$ and $v\in H^1(\mathcal T_h)$ then it holds:
	\begin{align}\label{3a}
	\tilde{a}(u,v-u)-\skalarehzero{u}{v-u}+j(v)-j(u)\geq(f,v-u).
	\end{align}
\end{Lemma}
{\bf Proof.} It follows from Corollary \ref{Corollary_1} that $u$ satisfies (\ref{1})-(\ref{3}). We multiply (\ref{1}) by $v-u$ and integrate over $\Omega$. Next, using Green formula for every element $K\in\Th$ separately, we obtain
\begin{align}\label{4}
\tilde{a}(u,v-u)-\sum_{K\in\Th}\int_{\partial K}\nabla u\cdot n(v-u)\,ds=\int_\Omega f(v-u)\,dx.
\end{align}  
On the other hand, applying formula (3.2) of \cite{Tutorial_DG}, and using the fact that the jump $[\nabla u]$ vanishes on internal edges, as $u\in H^2(\Omega)\subset C(\Omega)$, we obtain
\begin{align}\label{5}
&\sum_{K\in\Th}\int_{\partial K}\nabla u\cdot n(v-u)\,ds=\skalareh{u}{v-u}+\langle[\nabla u],\{v-u\}\rangle_{\mathcal E_h^0}\nonumber\\
&=\skalareh{u}{v-u}=\skalarehzero{u}{v-u}+\int_{\partial\Omega}\nabla u\cdot n(v-u)\,ds.
\end{align}
Applying (\ref{5}) in (\ref{4}), we get
\begin{align}\label{6}
\tilde{a}(u,v-u)-\skalarehzero{u}{v-u}
-\int_{\partial\Omega}\nabla u\cdot n(v-u)\,ds=\int_\Omega f(v-u)\,dx.
\end{align}
From (\ref{3}) and assumptions on $\lambda$ and $g$, we obtain the thesis. 
$\hfill{\Box}$ \\

Now we introduce notation of the virtual element method. For every $K\in \Th$ and fixed $k\in \nat_+$, we define
\begin{align*}
&B_k(\partial K)=\{v\in C^0(\partial K),\,\,\,v_{|e}\in{\mathbb P}_k(e)\,\,\,\,\,\mbox{for all edges} ~  e ~ \text{of}\,\,K \},\\[2mm]
&V_k(K)=\{v\in H^1(K):\,\,v_{|\partial K}\in B_k(\partial K)\,\,\text{and}\,\,\,\Delta v\in {\mathbb P}_k(K)  \}.
\end{align*} 
Here the symbol ${\mathbb P}_k$ stands for the set of polynomials of order $\leq k$. We also introduce the space of monomials of order $l\in\nat$ defined by
\begin{align*}
{\mathbb M}_{l}(K)=
\left\lbrace  
\left(\frac{x-x_K}{h_K}\right)^s, 
|s|\leq l\right\rbrace,\\
{\mathbb M}_{l}^*(K)=
\left\lbrace  
\left(\frac{x-x_K}{h_K}\right)^s, 
|s|= l\right\rbrace,
\end{align*}
where $x_K$ is a barycenter of and $s=(s_1,s_2)\in\nat^2$ is a multiindex. Next we define projection operator $\Pi_k^K\colon V_k(K)\to {\mathbb P}_k(K)$ as follows. For all $v\in V_k(K)$, $\Pi_k^Kv$ is the (unique!) element of ${\mathbb P}_k(K)$ that satisfies
\begin{align}
&(\nabla\Pi_k^Kv,\nabla p)_{0,K}=(\nabla v,\nabla p)_{0,K}\quad\mbox{for all} ~  p\in {\mathbb P}_k(K),\label{7}\\[2mm]
&\int_{\partial K}(v-\Pi_k^Kv)\,ds=0.\label{8}
\end{align}
Now we define the following subspace $W_k(K)$ of $V_k(K)$.
\begin{align*}
W_k(K)=\{ w\in V_k(K):\,\,(w-\Pi_k^Kw,q)_{0,K}=0,\,\,\,\,\,\,\mbox{for all} ~  q\in {\mathbb M}_{k-1}^*(K)\cup {\mathbb M}_{k}^*(K)\}
\end{align*}
and the space 
\begin{align*}
W_{DG}=\{v\in L^2(\Omega):\,\,\,v_{|K}\in W_k(K)\,\,\,\,\,\mbox{for all} ~  K\in \Th \}.
\end{align*}
We define the bilinear form $a_h^K\colon W_k(K)\times W_k(K)\to \real$ by 
\begin{align*}
a_h^K(u,v)=a^K(\Pi_k^Ku,\Pi_k^Kv)+S^K(u-\Pi_k^Ku,v-\Pi_k^Kv)\qquad\mbox{for all} ~  u,v\in W_k(K),
\end{align*}
where $S^K(u,v)$ is any symmetric, positive definite, bilinear form which satisfies:
\begin{align*}
c_* a^K(v,v)\leq S^K(v,v)\leq c^*a^K(v,v)\qquad\mbox{for all} ~  v\in W_k(K)\,\,\,\text{with}\,\,\,\Pi_k^Kv=0,
\end{align*}
where $c_*, c^*>0$ are constants independent of the parameter of discretization $h$.\\ 
We remark that $a_h^K$ satisfies the following properties
\begin{enumerate}
	\item $k$-consistency:
	\begin{align*}
a_h^K(p,v)=a^K(p,v)\qquad\mbox{for all} ~  p\in {\mathbb P}_k(K),\,\,\,v\in W_k(K),
	\end{align*}
	\item stability: 
	\begin{align}\label{8b}
	\alpha_* a^K(v,v)\leq a_h^K(v,v)\leq\alpha^* a^K(v,v)\qquad \mbox{for all} ~  v\in W_k(K),	
	\end{align}
	where $\alpha_*, \alpha^*>0$ are constants independent of $h$.
\end{enumerate}
We recall that the above notation is suitable for any $k\in\nat_+$. However, for the sake of simplicity, we restrict ourselves to the case $k=1$. Furthermore, we simplify notation and write $\Pi^K$ instead of $\Pi^K_1$ and $W(K)$ instead of $W_1(K)$. Hence, the space $W(K)$ can be characterized as follows
\begin{align*}
W(K)=\{w\in V_1(K):\,\,\, (w-\Pi^Kw,q)_{0,K}=0,\quad \mbox{for all} ~  q\in \mathbb{P}_1(K)\}.
\end{align*}  	
\begin{Corollary}\label{Corrolary_2}
For all $K\in\Th $ and for all $w\in W(K)$, the projection $\Pi^K w$ corresponds to the $L^2(K)$-projection of $w$ on the set $\mathbb{P}_1(K)$.
\end{Corollary}
Finally, we define the form 
\begin{align*}
\tilde{a}_h(u,v)=\sum_{K\in\Th}a_h^K(u,v)\qquad\mbox{for all} ~  v\in W_{DG}
\end{align*} 
and the approximation $f_h$ of the right hand side of Problem $(P)$
\begin{align}
\skalar{f_h,v_h}=\sum_{K\in \Th}\int_Kf\cdot\Pi^Kv_h\,dx\qquad\mbox{for all} ~  v_h\in H^1(\Th). \label{f_h}
\end{align} 
We introduce the following norms for $v\in H^2(\Th)$
\begin{align*}
&\|v\|^2_{2,DG}=\sum_{K\in\Th}\left(\|\nabla v\|^2_{0,K}+h_K^2|\nabla v|^2_{1,K}+\|v\|^2_{0,K} \right)+\sskalareh{v}{v},\\[2mm]
&\|v\|^2_{1,DG}=\sum_{K\in\Th}\left(\|\nabla v\|^2_{0,K}+\|v\|^2_{0,K} \right)+\sskalareh{v}{v}.
\end{align*}
We remark, that for a function $v\in H^2(\Th)$ which is piecewise polynomial, we have:
\begin{align*}
\|v\|_{2,DG}\simeq\|v\|_{1,DG}.
\end{align*} 
Now we define the bilinear form $B_h\colon W_{DG}\times W_{DG}\to \real$ by
\begin{align*}
B_h(v,w)=\tilde{a}_h(v,w)-\skalarehzero{\Pi v}{w}-\delta\skalarehzero{\Pi w}{v}+\gamma\sskalareh{v}{w},
\end{align*}
where $\delta\in \{-1,0,1\}$ and $\gamma$ is a constant whose value will be fixed letter. The symbol $\Pi$ in the above formula is understood as follows
\begin{align*}
\Pi v_{|K}=\Pi^K v\quad \text{for all}\,\,K\in\Th, v\in W_{DG}.
\end{align*}
The following lemma provides a kind of coercivity and continuity of the form $B_h$.
\begin{Lemma}\label{Lemma_2}
	The form $B_h$ satisfies
	\begin{align}
	& B_h(v,v)\geq M_s\|v\|^2_{1,DG}\qquad\mbox{for all} ~  v\in W_{DG}\label{10}
	\end{align}
	and we have
	\begin{align}
	&|\tilde{a}_h(w,v)|\leq c\|v\|_{1,DG}\|w\|_{1,DG}\qquad\mbox{for all} ~  v,w\in W_{DG}\label{11a},\\[2mm]
	&|\sskalareh{v}{w}|\leq c\|v\|_{1,DG}\|w\|_{1,DG}\qquad\mbox{for all} ~  v,w\in W_{DG}\label{11b},\\[2mm]
	&\skalareh{v}{w}\leq c \|v\|_{2,DG}\|w\|_{1,DG}\qquad\mbox{for all} ~  v\in H^2(\Th), w\in H^1(\Th)\label{12}.
	\end{align}
Moreover if $v$ is piecewise polynomial, then 
\begin{align}\label{13}
\skalareh{v}{w}\leq c \|v\|_{1,DG}\|w\|_{1,DG}\qquad\mbox{for all} ~  w\in H^1(\Th).
\end{align}	
\end{Lemma}

\noindent To proceed further, we assume that all virtual elements are affine equivalent, see  \cite[Definition 3.4.1]{BRESCO}, which is possible due to assumption $(C_0)$. 

Before the proof of Lemma \ref{Lemma_2}, we recall the scaled trace inequality  \cite[(10.3.8)]{BRESCO}
\begin{align}\label{14}
\|v\|^2_{0,\partial K}\leq c(h_K^{-1}\|v\|^2_{0,K}+h_K|v|^2_{1,K})\qquad\mbox{for all} ~  v\in H^1(K)
\end{align}
and the inverse inequality \cite[Lemma 4.5.3]{BRESCO}
\begin{align}\label{14a}
|v|^2_{1,K}\leq c h_K^{-2}\|v\|^2_{0,K}\qquad\mbox{for all} ~  v\in{\mathbb P}_1(K).     
\end{align}
Now we pass to the proof of Lemma \ref{Lemma_2}.\\

\noindent {\bf Proof}. First, we prove (\ref{12}). To this end, we have for $v\in H^2(\Th), w\in H^1(\Th)$
\begin{align}\label{15}
&\skalareh{v}{w}=\sum_{e\in{\mathcal E_h}}\int_e\{\nabla v\}[w]ds\nonumber\\[2mm]
&\leq\left(h\sum_{e\in{\mathcal E_h}}\int_e\{\nabla v\}^2\,ds\right)^{\frac{1}{2}}\left(\frac{1}{h}\sum_{e\in{\mathcal E_h}}\int_e [w]^2\,ds\right)^{\frac{1}{2}}.
\end{align}
Moreover, we estimate
\begin{align}\label{16}
&\sum_{e\in{\mathcal E_h}}\int_e\{\nabla v\}^2\,ds=\sum_{e\in{\mathcal E_h^\partial}}\int_e|\nabla v|^2\,ds+\frac{1}{4}\sum_{e\in{\mathcal E_h^0}}\int_e|\nabla v^++\nabla v^-|^2\,ds\nonumber\\[2mm]
&\leq\sum_{e\in{\mathcal E_h^\partial}}\int_e|\nabla v|^2\,ds+\frac{1}{2}\sum_{e\in{\mathcal E_h^0}}\int_e\left(|\nabla v^+|^2+|\nabla v^-|^2\right)\,ds\nonumber\\[2mm]
&\leq c\sum_{K\in\Th}\int_{\partial K}|\nabla v|^2\,ds=c\sum_{K\in\Th}\|\nabla v\|^2_{0,\partial K}.
\end{align}
As $v\in H^2(\Th)$, we $|\nabla v|\in H^1(\Th)$. Hence, we can apply (\ref{14}) for $|\nabla v|$ in (\ref{16}) to obtain
\begin{align}\label{17}
\left(h\sum_{e\in{\mathcal E_h}}\int_e\{\nabla v\}^2\,ds\right)^{\frac{1}{2}}\leq c\left(\sum_{K\in\Th}\|\nabla v\|^2_{0,K}+h_K^2\sum_{K\in\Th}|\nabla v|^2_{1,K}\right)^{\frac{1}{2}}\leq c\|v\|_{2,DG}.
\end{align}
On the other hand, we have
\begin{align}\label{18}
\left(\frac{1}{h}\sum_{e\in{\mathcal E_h}}\int_e [w]^2\,ds\right)^{\frac{1}{2}}\leq \left(\sum_{e\in{\mathcal E_h}}\frac{1}{|e|}\int_e[w]^2\,ds\right)^{\frac{1}{2}}=\sskalareh{w}{w}^{\frac{1}{2}}\leq c\|w\|_{1,DG}.
\end{align}
Now, from (\ref{15}), (\ref{17}) and (\ref{18}), we obtain (\ref{12}).

Now we prove (\ref{11a}). Recall, that the form $a_h^K$ is bilinear, symmetric and positive definite. Hence, using also (\ref{8b}), we have
\begin{align}\label{19}
&\tilde{a}_h(v,w)=\sum_{K\in\Th}a_h^K(v,w)\leq\left(\sum_{K\in\Th}a_h^K(v,v)\right)^\frac{1}{2}\left(\sum_{K\in\Th}a_h^K(w,w)\right)^\frac{1}{2}\nonumber\\[2mm]
&\leq\alpha^*\left(\sum_{K\in\Th}a^K(v,v)\right)^\frac{1}{2}\left(\sum_{K\in\Th}a^K(w,w)\right)^\frac{1}{2}\nonumber\\[2mm]
&=\alpha^*\left(\sum_{K\in\Th}(\|\nabla v\|^2_{0,K}+\|v\|^2_{0,K})\right)^\frac{1}{2}\left(\sum_{K\in\Th}\|\nabla w\|^2_{0,K}+\|w\|^2_{0,K}\right)^\frac{1}{2}\nonumber\\[2mm]
&\leq\alpha^*\|v\|_{1,DG}\|w\|_{1,DG}.
\end{align}   
Next, we estimate
\begin{align}\label{20}
&|\sskalareh{v}{w}|=\sum_{e\in{\mathcal E}_h}\frac{1}{|e|}\int_e|[v][w]|ds\nonumber\\[2mm]
&\leq\left(\sum_{e\in{\mathcal E}_h}\frac{1}{|e|}\int_e|[v]|^2ds\right)^{\frac{1}{2}}\left(\sum_{e\in{\mathcal E}_h}\frac{1}{|e|}\int_e|[w]|^2ds\right)^{\frac{1}{2}}\leq\|v\|_{1,DG}\|w\|_{1,DG},
\end{align}
which shows (\ref{11b}). Now we pass to the proof of (\ref{10}). From (\ref{8b}), we get
\begin{align}\label{21}
\tilde{a}_h(v,v)=\sum_{K\in\Th}a_h^K(v,v)\geq \alpha_*\sum_{K\in\Th}a^K(v,v)=\alpha_*\sum_{K\in\Th}\left(\|\nabla v\|^2_{0,K}+\|v\|^2_{0,K}\right). 
\end{align} 
Next, from (\ref{15}), we have
\begin{align}\label{22}
\skalarehzero{\Pi v}{v}\leq\left(h\sum_{e\in{\mathcal E}_h^0}\int_e\{\nabla\Pi v\}^2\,ds\right)^{\frac{1}{2}}\left(\frac{1}{h}\sum_{e\in{\mathcal E}_h^0}\int_e[v]^2\,ds\right)^{\frac{1}{2}}.
\end{align}
Proceeding in a similar way as in the proof of (\ref{17}) and using (\ref{14a}), we get
\begin{align}\label{23}
&\left(h\sum_{e\in{\mathcal E_h}}\int_e\{\nabla\Pi v\}^2\,ds\right)^{\frac{1}{2}}\leq c\left(\sum_{K\in\Th}\|\nabla\Pi v\|^2_{0,K}+h_K^2\sum_{K\in\Th}|\nabla\Pi v|^2_{1,K}\right)^{\frac{1}{2}}\nonumber\\[2mm]
&\leq c\left(\sum_{K\in\Th}\|\nabla\Pi v\|^2_{0,K}\right)^{\frac{1}{2}}.
\end{align}
Now we estimate the norm $\|\nabla\Pi v\|_{0,K}$. By the definition (\ref{7}), we have
\begin{align*}
\|\nabla\Pi v\|_{0,K}^2=(\nabla\Pi v,\nabla\Pi v)_{0,K}=(\nabla\Pi v,\nabla v)_{0,K}\leq\|\nabla\Pi v\|_{0,K}\|\nabla v\|_{0,K}.
\end{align*}
Hence
\begin{align}\label{24}
\|\nabla\Pi v\|_{0,K}\leq\|\nabla v\|_{0,K}.
\end{align}
From (\ref{23}) and (\ref{24}), we have
\begin{align}\label{25}
&\left(h\sum_{e\in{\mathcal E_h}}\int_e\{\nabla\Pi v\}^2\,ds\right)^{\frac{1}{2}}\leq c\left(\sum_{K\in\Th}\|\nabla v\|^2_{0,K}\,ds\right)^{\frac{1}{2}}\nonumber\\[2mm]
&\leq c\left(\sum_{K\in\Th}\|\nabla v\|^2_{0,K}+\|v\|^2_{0,K}\,ds\right)^{\frac{1}{2}}.
\end{align}
From (\ref{22}) and (\ref{25}), we have
\begin{align}\label{26}
&-(1+\delta)\skalarehzero{\Pi v}{v}\nonumber\\[2mm]
&\geq-(1+\delta)\tilde{c}\left(\sum_{K\in\Th}\|\nabla v\|^2_{0,K}+\|v\|^2_{0,K}\,ds\right)^{\frac{1}{2}}\left(\sskalarehzero{v}{v}\right)^{\frac{1}{2}}
\end{align}
with $\tilde{c}>0$. From (\ref{21}), (\ref{26}) and definition of $B_h$, we get
\begin{align}\label{27}
B_h(v,v)\geq\alpha_*x^2+\gamma y^2-(1+\delta)\tilde{c}xy,
\end{align}
where
\begin{align*}
x:=\left(\sum_{K\in\Th}\|\nabla v\|^2_{0,K}+\|v\|^2_{0,K}\,ds\right)^{\frac{1}{2}}, \quad y:=\left(\sskalareh{v}{v}\right)^{\frac{1}{2}}.
\end{align*}
We remark that for $\gamma>\frac{(1+\delta)^2\tilde{c}^2}{\alpha_*}$, there exists $M_s>0$ such that
\begin{align*}
\alpha_*x^2+\gamma y^2-(1+\delta)\tilde{c}xy\geq M_s(x^2+y^2)=M_s\|v\|^2_{1,DG}.
\end{align*}
Combining the latter with (\ref{27}), we conclude (\ref{10}). From (\ref{12}) and (\ref{14a}), we obtain (\ref{13}), which completes the proof. $\hfill{\Box}$\\

\section{Approximation and convergence}

\noindent In this section we deal with the discrete approximation of Problem $(V)$. Consider the following problem.\\

\noindent {\bf Problem ($V_{DG}$).} {\it Find $u_h\in W_{DG}$ such that
	\begin{align}\label{28}
	B_h(u_h,v_h-u_h)+j(v_h)-j(u_h)\geq \skalar{f_h,v_h-u_h}\qquad \mbox{for all} ~ \,\, v_h\in W_{DG}.
	\end{align}
} 
We start with an existence and uniqueness result for Problem ($V_{DG}$). 
\begin{Proposition}\label{Proposition_1}
Problem ($V_{DG}$) has a unique solution.	
\end{Proposition}
{\bf Proof.} It follows from Lemma \ref{Lemma_2} that $B_h$ is $V_{DG}$-elliptic and continuous with respect to the norm $\|\cdot\|_{1,DG}$. Moreover, the functional $W_{DG}\ni v_h\to \skalar{f_h,v_h}\in\real$ is continuous and $j$ is l.s.c. Now the thesis follows from Theorem \ref{Theorem_0}.
$\hfill{\Box}$\\

The main goal of this section is to derive a linear bound for an error between the solution of Problem $(V)$ and the one of Problem  ($V_{DG}$). The final result, Theorem \ref{Theorem_2} will be preceded by Lemma \ref{Lemma_4} and Theorem \ref{Theorem_1}. 

\begin{Lemma}\label{Lemma_4}
	For every $w\in W_{DG}$, we have
	\begin{align}\label{A1}
	\|\Pi w\|_{1,DG}\leq c\|w\|_{1,DG}.
	\end{align}
\end{Lemma}

\noindent{\bf Proof.} Let $w\in W_{DG}$. Then, for every $K\in\Th$, we have $w_{|K}\in W(K)$. Hence, from Corollary \ref{Corrolary_2}, we have
\begin{align}\label{eq_2}
\|\Pi^Kw\|_{0,K}\leq\|w\|_{0,K}.
\end{align} 

\noindent 
From (\ref{24}) and (\ref{eq_2}), we have
\begin{align}\label{A2}
\sum_{K\in\Th}\|\nabla\Pi^Kw\|^2_{0,K}+\sum_{K\in\Th}\|\Pi^Kw\|^2_{0,K}\leq\sum_{K\in\Th}\|\nabla w\|^2_{0,K}+\sum_{K\in\Th}\|w\|^2_{0,K}.
\end{align}
Moreover, we obtain
\begin{align}\label{A3}
&\sskalareh{\Pi w}{\Pi w}=\sum_{e\in{\mathcal E_h}}\frac{1}{|e|}\int_e|[\Pi w]^2|\,ds=\sum_{e\in{\mathcal E_h}}\frac{1}{|e|}\int_e\left([\Pi w-w]+[w]\right)^2\,ds\nonumber\\[2mm]
&\leq 2\sum_{e\in{\mathcal E_h}}\frac{1}{|e|}\int_e[\Pi w-w]^2\,ds+2\sskalareh{w}{w}.
\end{align}
Using (\ref{14}), the Poincare inequality and (\ref{24}), we get
\begin{align}\label{A4}
&2\sum_{e\in{\mathcal E_h}}\frac{1}{|e|}\int_e[\Pi w-w]^2\,ds\nonumber\\[2mm]
&\leq c\frac{1}{e_K^{min}}\sum_{K\in\Th}\int_{\partial K}|\Pi w-w|^2\,ds\leq c\frac{1}{\gamma_2h_K}\sum_{K\in\Th}\|\Pi w-w\|^2_{0,\partial K}\nonumber\\[2mm]
&\leq\frac{c}{\gamma_2}\sum_{K\in\Th}\left(\frac{1}{h_K^2}\|\Pi w-w\|^2_{0,K}+\|\nabla\Pi w-\nabla w\|^2_{0,K}\right)\nonumber\\[2mm]
&\leq c\sum_{K\in\Th}\|\nabla\Pi w-\nabla w\|^2_{0,K}\leq 2c\sum_{K\in\Th}\left(\|\nabla\Pi w\|^2_{0,K}+\|\nabla w\|^2_{0,K}\right)\nonumber\\[2mm]
&\leq 4c\sum_{K\in\Th}\|\nabla w\|^2_{0,K}.
\end{align}
Here $e_K^{min}$ denotes the length of the shortest edge of element $K$ and $\gamma_2$ is the constant from assumption $(C_2)$.  Note, that due to (\ref{8}) and the definition of the space $W_{DG}$, we could use Poincare-Friedrich's inequality(e.g., \cite[Theorem 1.1]{Farwig}), since average of $\Pi w-w$ on $K$ is zero.
From (\ref{A2})-(\ref{A4}), we get (\ref{A1}), which completes the proof. 
$\hfill{\Box}$

\begin{Theorem}\label{Theorem_1}
	Let $u\in H^2(\Omega)$ be a solution of Problem $(V)$ and $u_h$ be a solution of Problem $(V_{DG})$. Then, for every approximation $u_I$ of $u$ in $W_{DG}$ and for every approximation $u_{\pi}$ of $u$ that is piecewise in ${\mathbb P_1}$, we have
	\begin{align}\label{30}
	&\|u-u_h\|^2_{1,DG}\leq\nonumber\\
	 &c(\|u-u_I\|^2_{1,DG}+\|u-u_\pi\|^2_{2,DG}+\|f-f_h\|^2_{W^*_{1,DG}} +|R(u,u_I)|),
	\end{align} 
where the residual term $R$ is defined by
\begin{align}\label{residuum}
&R(u,v)=\tilde{a}(u,v-u)-\skalarehzero{u}{v-u}\nonumber\\
&\qquad\qquad\qquad+j(v)-j(u)-(f,v-u)\qquad\mbox{for all} ~  u,v\in W_{DG}.
\end{align}	
\end{Theorem}
{\bf Proof.} Let $u_I\in W_{DG}$ be given and $\eta_h=u_h-u_I$. Using (\ref{10}), we have
\begin{align}\label{31}
M_s\|\eta_h\|^2_{1,DG}\leq B_h(\eta_h,\eta_h)=B_h(u_h,\eta_h)-B_h(u_I,\eta_h).
\end{align}
It follows from (\ref{28}) that (taking $v_h=u_I$)
\begin{align}\label{32}
B_h(u_h,\eta_h)=B_h(u_h,u_h-u_I)\leq j(u_I)-j(u_h)+\skalar{f_h,u_h-u_I}.
\end{align}
From (\ref{31}), (\ref{32}) and by the definition of $B_h$, we get
\begin{align}\label{33}
M_s\|\eta_h\|^2_{1,DG}&\leq j(u_I)-j(u_h)+\skalar{f_h,u_h-u_I}-\tilde{a}_h(u_I,\eta_h)\nonumber\\[1mm]
&+\skalarehzero{\Pi u_I}{\eta_h}+\delta\skalarehzero{\Pi \eta_h }{u_I}-\gamma\sskalareh{u_I}{\eta_h}\nonumber\\[1mm]
&=j(u_I)-j(u_h)+\skalar{f_h,u_h-u_I}-\tilde{a}_h(u_I,\eta_h)\nonumber\\[1mm]
&+\skalarehzero{\Pi (u_I-u_\pi)}{\eta_h}+\skalarehzero{\Pi u_\pi}{\eta_h}\nonumber\\[1mm]
&+\delta\skalarehzero{\Pi \eta_h }{u_I}-\gamma\sskalareh{u_I}{\eta_h}.
\end{align}
Moreover, using the consistency property of $a_h$, we have
\begin{align}\label{34}
&\tilde{a}_h(u_I,\eta_h)=\sum_{K\in\Th}a_h^K(u_I,\eta_h)=\sum_{K\in\Th}\left(a_h^K(u_I-u_\pi,\eta_h)+a_h^K(u_\pi,\eta_h)\right)\nonumber\\[2mm]
&=\sum_{K\in\Th}a_h^K(u_I-u_\pi,\eta_h)+\sum_{K\in\Th}a^K(u_\pi,\eta_h)\nonumber\\[2mm]
&=\sum_{K\in\Th}a_h^K(u_I-u_\pi,\eta_h)+\sum_{K\in\Th}a^K(u_\pi-u,\eta_h)+\sum_{K\in\Th}a^K(u,\eta_h)\nonumber\\[2mm]
&=\sum_{K\in\Th}a_h^K(u_I-u_\pi,\eta_h)+\sum_{K\in\Th}a^K(u_\pi-u,\eta_h)+\tilde{a}(u,\eta_h).
\end{align}
From (\ref{33}) and (\ref{34}), we get
\begin{align}\label{35}
M_s\|\eta_h\|^2_{1,DG}&=j(u_I)-j(u_h)+\skalar{f_h,u_h-u_I}-\sum_{K\in\Th}a_h^K(u_I-u_\pi,\eta_h)\nonumber\\[1mm]
&-\sum_{K\in\Th}a^K(u_\pi-u,\eta_h)-\tilde{a}(u,\eta_h)+\skalarehzero{\Pi (u_I-u_\pi)}{\eta_h}\nonumber\\[1mm]
&+\skalarehzero{\Pi u_\pi}{\eta_h}+\delta\skalarehzero{\Pi \eta_h }{u_I}-\gamma\sskalareh{u_I}{\eta_h}.
\end{align}
Since $u$ is a solution of Problem $(V)$ and $u\in H^2(\Omega)$, it follows from Lemma \ref{Lemma_1} that (\ref{3a}) holds. From (\ref{3a}), we have (for $v=u_h$)
\begin{align}\label{36}
-\tilde{a}(u,\eta_h)=-\tilde{a}(u,u_h-u_I)=-\tilde{a}(u,u_h-u)-\tilde{a}(u,u-u_I)\leq\nonumber\\[2mm]
-\skalarehzero{u}{u_h-u}+j(u_h)-j(u)+\skalar{f,u-u_h}-\tilde{a}(u,u-u_I)
\end{align}
Applying (\ref{36}) and the definition of the residual term (\ref{residuum}) in (\ref{35}), we get
\begin{align}\label{38}
M_s\|\eta_h\|^2_{1,DG}&\leq R(u,u_I)-\skalarehzero{u}{u_h-u_I} +\skalar{f,u_I-u_h}+\skalar{f_h,u_h-u_I}\nonumber\\[1mm]
&-\sum_{K\in\Th}a_h^K(u_I-u_\pi,\eta_h)-\sum_{K\in\Th}a^K(u_\pi-u,\eta_h)\nonumber\\[1mm]
&+\skalarehzero{\Pi (u_I-u_\pi)}{\eta_h}+\skalarehzero{\Pi u_\pi}{\eta_h}\nonumber\\[1mm]
&+\delta\skalarehzero{\Pi \eta_h }{u_I}-\gamma\sskalareh{u_I}{\eta_h}=\nonumber\\[1mm]
&R(u,u_I)+\skalarehzero{(u_\pi-u)}{\eta_h}+\skalarehzero{\Pi (u_I-u_\pi)}{\eta_h}\nonumber\\[1mm]
&+\delta\skalarehzero{\Pi \eta_h }{u_I}-\gamma\sskalareh{u_I}{\eta_h}+\skalar{f_h-f,\eta_h}\nonumber\\[1mm]
&-\sum_{K\in\Th}a_h^K(u_I-u_\pi,\eta_h)-\sum_{K\in\Th}a^K(u_\pi-u,\eta_h).
\end{align}
In the last step, we have used the fact $\Pi u_\pi=u_\pi$.
Now we estimate each term of the right hand side of (\ref{38}). From (\ref{12}), we have
\begin{align}\label{39}
\skalarehzero{ (u_\pi-u)}{\eta_h}\leq c\|u_\pi-u\|_{2,DG}\|\eta_h\|_{1,DG}.
\end{align}
From (\ref{13}) and Lemma \ref{Lemma_4}, we deduce
\begin{align}\label{40}
\skalarehzero{\Pi (u_I-u_\pi)}{\eta_h}&\leq c\|\Pi(u_I-u_\pi)\|_{1,DG}\|\eta_h\|_{1,DG}\nonumber\\[1mm]
&\leq c\|u_I-u_\pi\|_{1,DG}\|\eta_h\|_{1,DG}.
\end{align}
Since $u\in H^2(\Omega)$, it has no jump, so $[u]=0$. Hence, from (\ref{13}) and Lemma \ref{Lemma_4}, we have
\begin{align}\label{41}
\skalarehzero{\Pi \eta_h }{u_I}&=\skalarehzero{\Pi \eta_h }{u_I-u}\leq c\|\Pi \eta_h\|_{1,DG}\|u_I-u\|_{1,DG}\nonumber\\[1mm]
&\leq c\|\eta_h\|_{1,DG}\|u_I-u\|_{1,DG}.
\end{align} 
Similarly from (\ref{11b}), we obtain
\begin{align}\label{42}
&-\sskalareh{u_I}{\eta_h}=-\sskalareh{u_I-u}{\eta_h}\leq |\sskalareh{u_I-u}{\eta_h}|\nonumber\\[1mm]
&\leq c\|u_I-u\|_{1,DG}\|\eta_h\|_{1,DG}.
\end{align}
Furthermore from \eqref{f_h}, we get
\begin{align}\label{43}
\skalar{f_h-f,\eta_h}\leq\|f_h-f\|_{W_{1,DG}^*}\|\eta_h\|_{1,DG}
\end{align}
The condition (\ref{11a}) implies
\begin{align}\label{44}
&-\sum_{K\in\Th}a_h^K(u_I-u_\pi,\eta_h)=-\tilde{a}_h(u_I-u_\pi,\eta_h)\nonumber\\
&\leq|\tilde{a}_h(u_I-u_\pi,\eta_h)|\leq c\|u_I-u_\pi\|_{1,DG}\|\eta_h\|_{1,DG},
\end{align}
which entails
\begin{align}\label{45}
-\sum_{K\in\Th}a^K(u_\pi-u,\eta_h)\leq c\|u_\pi-u\|_{1,DG}\|\eta_h\|_{1,DG}.
\end{align}
From (\ref{38})-(\ref{45}), we have
\begin{align}\label{46}
&\|\eta_h\|^2_{1,DG}\leq|R(u,u_I)|+\nonumber\\[2mm]
&c(\|u_\pi-u\|_{2,DG}+\|u_I-u_\pi\|_{1,DG}+\|u_I-u\|_{1,DG}+\|f_h-f\|_{W_{1,DG}^*})\|\eta_h\|_{1,DG}.
\end{align}
From (\ref{46}), we easily get
\begin{align}\label{47}
&\|\eta_h\|^2_{1,DG}\leq\nonumber\\[2mm]
&c(|R(u,u_I)|+\|u_\pi-u\|^2_{2,DG}+\|u_I-u_\pi\|^2_{1,DG}+\|u_I-u\|^2_{1,DG}+\|f_h-f\|^2_{W_{1,DG}^*}).
\end{align}
From the triangle inequality, we easily deduce
\begin{align}
&\|u-u_h\|^2_{1,DG}\leq 2(\|u-u_I\|^2_{1,DG}+\|\eta_h\|^2_{1,DG}),\label{48}\\[2mm]
&\|u_I-u_\pi\|^2_{1,DG}\leq 2(\|u_I-u\|^2_{1,DG}+\|u-u_\pi\|^2_{1,DG}).\label{49}
\end{align}
Applying (\ref{48}) and (\ref{49}) in (\ref{47}), we obtain (\ref{30}), which completes the proof. $\hfill{\Box}$\\

Now we are in a position to present the main result of this paper. 
\begin{Theorem}\label{Theorem_2}
Let $u\in H^2(\Omega)$ be the solution of Problem ($V$), such that $u_{|e}\in H^2(e)$ $\mbox{for all} ~  e\in {\mathcal E}^{\partial}$. Let $u_h$ be the solution of Problem $(V_{DG})$. Then, under the assumptions $(C_0)-(C_2)$ it holds 
\begin{align}\label{49A}
\|u-u_h\|_{1,DG}\leq ch,
\end{align}
where $c$ is a positive constant independent of $h$.
\end{Theorem}
{\bf Proof}. We use Theorem \ref{Theorem_1} and we estimate the right hand side of (\ref{30}). In particular, we take in (\ref{30}) the element $u_I$ defined as follows.
\begin{align*}
u_{I|K}=I_Ku\quad\mbox{for} ~  K\in\Th,
\end{align*}
where $I_Ku\in W_{DG}$ is such that $(I_Ku)(x_i)=u(x_i)$ for all vertices $x_i$ of $K$. Hence $u_I$ is piecewise linear on every $e\in{\mathcal E}$ and 
\begin{align}\label{50}
[u_I]=0\quad\mbox{for all} ~  e\in {\mathcal E^0_h}.
\end{align}
Moreover, since $u\in H^2(\Omega)$, we have 	
\begin{align}\label{51}
[u]=0\quad\mbox{for all} ~  e\in {\mathcal E^0_h}.
\end{align}
From (\ref{50})-(\ref{51}), we get	
\begin{align}\label{52}
[u-u_I]=0\quad\mbox{for all} ~  e\in {\mathcal E^0_h}.
\end{align}
From (\ref{52}), we obtain
\begin{align}\label{53}
&\skalar{[u-u_I],[u-u_I]}_{\mathcal E_h}=\sum_{e\in{\mathcal E}^\partial}\frac{1}{|e|}\int_e|u-u_I|^2ds\nonumber\\[2mm]
&=\sum_{e\in{\mathcal E}^\partial_h}\frac{1}{|e|}\|u-u_I\|^2_{L^2(e)}\leq c\sum_{e\in{\mathcal E}^\partial_h}\frac{1}{|e|}|e|^4|u|^2_{H^2(e)}\nonumber\\[2mm]
&\leq h^3|u|^2_{H^2(\partial\Omega)}.
\end{align} 	
Now we use the following property of the interpolation operator $I_K\colon H^{k+1}(K)\to W_k(K)$. Under the assumptions $(C_0)-(C_2)$, we have, see \cite[Theorem 4.4.4]{BRESCO}
\begin{align}\label{54}
\|v-I_Kv\|_{0,K}+h_K|v-I_Kv|_{1,K}\leq c h_K^{k+1}|v|_{k+1,K}.
\end{align}
We underline that the constant $c$ in (\ref{54}) does not depend on element $K$. It follows from the affine equivalence of all elements and from \cite[Proposition 4.4.11]{BRESCO}.
From (\ref{54}) we estimate
\begin{align}
&\|u-u_I\|^2_{0,K}\leq ch^4|u|^2_{2,K},\label{55}\\[2mm]
&|u-u_I|^2_{1,K}\leq ch^2|u|^2_{2,K}.\label{56}
\end{align}
From (\ref{53}), (\ref{55}) and (\ref{56}), we get
\begin{align}\label{57}
\|u-u_I\|^2_{1,DG}\leq c(h^3+h^4+h^2).
\end{align}
We define the $L^2$-projection $\Pi^0\colon W_1(K)\to{\mathbb P}_1(K)$ (it can be extended to $\Pi^0\colon H^1(K)\to{\mathbb P}_1(K)$), by
\begin{align*}
(\Pi^0v,m)_K=(\Pi v,m)\quad\mbox{for all} ~  m\in{\mathbb P}_1(K).
\end{align*}
Now we take $u_\pi:=\Pi^0 u$ in (\ref{30}) and we use the following result (\cite[Lemma 4.3.8]{BRESCO}:
\begin{align}\label{58}
\|u-u_\pi\|_{0,K}+h_K|u-u_\pi|_{1,K}\leq ch_K^2\|u\|_{H^2(K)}.
\end{align}
Using (\ref{58}), we estimate the term $\|u-u_\pi\|^2_{2,DG}$. Namely
\begin{align}
&\sum_{K\in\Th}\|u-u_\pi\|^2_{0,K}\leq ch^4\|u\|^2_{H^2(\Omega)},\label{59}\\[2mm]
&\sum_{K\in\Th}\|\nabla(u-u_\pi)\|^2_{0,K}\leq ch^2\|u\|^2_{H^2(\Omega)},\label{60}\\[2mm]
&\sum_{K\in\Th}h_K^2|\nabla(u-u_\pi)|^2_{1,K}=\sum_{K\in\Th}h_K^2|\nabla u|^2_{1,K}\leq h^2\|u\|^2_{H^2(\Omega)}\label{61}.
\end{align}
In the last estimate we use the fact that $u_\pi\in {\mathbb P}_1$ so all its second order derivatives vanish. Similarly as in (\ref{A4}) in the proof of Lemma \ref{Lemma_4}, we obtain (using (\ref{59}) and (\ref{60}))
\begin{align}\label{62}
&\skalar{[u-u_\pi][u-u_\pi]}_{\mathcal E_h}\leq c\sum_{K\in\Th}\left(\frac{1}{h^2}\|u-u_\pi\|^2_{0,K}+|u-u_\pi|^2_{1,K} \right)\nonumber\\[2mm]
&\leq ch^2\|u\|^2_{H^2(\Omega)}.
\end{align}
 From (\ref{59})-(\ref{62}), we get
 \begin{align}\label{62A}
 \|u-u_\pi\|^2_{2,DG}\leq ch^2\|u\|^2_{H^2(\Omega)}.
 \end{align}
It follows from (\ref{f_h}) that
\begin{align}\label{63}
\|f_h-f\|^2_{W_{1,DG}^*}\leq ch^2\|f\|^2_{L^2(\Omega)}.
\end{align}
It remains to estimate the residual term. Using (\ref{6}) and the definition of $R$, we see that
\begin{align}\label{64}
&|R(u,u_I)|=\left|\int_{\partial \Omega}\nabla u\cdot n(u_I-u)ds+j(u_I)-j(u)\right|\nonumber\\[2mm]
&=\left|\int_{\Gamma_2}-g\lambda(u_I-u)ds+\int_{\Gamma_2}g|u_I|-\int_{\Gamma_2}g|u|ds\right|\nonumber\\[2mm]
&=\left|\int_{\Gamma_2}g(|u_I|-\lambda u_I)ds\right|=\left|\int_{\Gamma_2}g(|u_I|-|u|+\lambda u-\lambda u_I)ds\right|\nonumber\\[2mm]
&\leq\int_{\Gamma_2}g\left|\,|u_I|-|u|+\lambda u-\lambda u_I\,\right|\,ds\nonumber\\[2mm]
&\leq\int_{\Gamma_2}g(|u_I-u|+|\lambda| |u-u_I|)ds\leq\int_{\Gamma_2}g(2|u_I-u|)ds\nonumber\\[2mm]
&\leq 2g\int_{\Gamma_2}|u_I-u|ds=2g\|u_I-u\|_{L^1(\Gamma_C)}.
\end{align} 
By the standard interpolation properties, we have
\begin{align*}
\|u_I-u\|_{L^1(\Gamma_C)}\leq ch^2|u|_{H^2(\Gamma_2)}.
\end{align*}
Combining it with (\ref{64}), we get
\begin{align}\label{65}
|R(u,u_I)|\leq ch^2.
\end{align}
From Theorem \ref{Theorem_1}, (\ref{57}), (\ref{62A}), (\ref{63}) and (\ref{65}), we obtain (\ref{49A}), which completes the proof of the theorem.
$\hfill{\Box}$ 

\section{Acknowledgements}\label{acknowledgement}
The research was
supported by the European
Union's Horizon 2020 Research and Innovation Programme under
the Marie Sklodowska-Curie grant agreement No. 823731 CONMECH,
the Ministry of Science and Higher Education of Republic of Poland
under Grant Nos. 4004/GGPJII/H2020/2018/0 and 440328/PnH2/2019,
and the National Science Centre of Poland under Project No. 2021/41/B/ST1/01636.

\end{document}